\documentclass[a4paper,12pt]{article}
\usepackage[T1]{fontenc}
\usepackage{amsmath,amssymb,amsthm}
\usepackage{microtype}
\usepackage{hyperref}
\hypersetup{hidelinks}
\allowdisplaybreaks[1]
\newtheorem{theorem}{Theorem}
\newtheorem{lemma}[theorem]{Lemma}
\newtheorem{proposition}[theorem]{Proposition}
\newtheorem{corollary}[theorem]{Corollary}
\newtheorem{mainthm}{Theorem}

\theoremstyle{definition}
\newtheorem{definition}[theorem]{Definition}
\theoremstyle{remark}

\DeclareMathOperator{\Aut}{Aut}
\DeclareMathOperator{\Sym}{Sym}
\DeclareMathOperator{\Syl}{Syl}
\DeclareMathOperator{\Ker}{Ker}

\DeclareMathOperator{\IYB}{IYB}
\DeclareMathOperator{\id}{id}
\newcommand{\Z}{\mathrm Z}
\newcommand{\C}{\mathrm C}
\newcommand{\N}{\mathrm N}

\title{Finite soluble groups with Sylow subgroups of class at most two are $\IYB$-groups}

\author{H. Meng\thanks{E-mail: hymeng2009@shu.edu.cn. Department of Mathematics
		and Newtouch Center for Mathematics of Shanghai University,
		Shanghai 200444, P. R. China. },~A. Ballester-Bolinches\thanks{E-mail: Adolfo.Ballester@uv.es. Departament de Matem\'{a}tiques, Universitat de Val\'{e}ncia, Dr. Moliner 50, 46100 Burjassot, Val\'{e}ncia,
		Spain.},~P. Jim\'{e}nez-Seral\thanks{E-mail: paz@unizar.es. Departamento de Matem\'{a}ticas, Universidad de Zaragoza, Pedro Cerbuna, 12, 50009 Zaragoza, Spain.}
}
\date{}

\begin{document}
\maketitle

\begin{abstract}
We prove that every finite soluble group whose Sylow subgroups have
nilpotency class at most two is an involutive Yang--Baxter group ($\IYB$-group).
This answers a question of Ced\'o and Okni\'nski \cite{CedoOkninski2025}.
The main method of the $\IYB$-structures has been reflected in \cite{MengBallesterJimenez2026}.


\end{abstract}

\section{Introduction}

A set-theoretic solution of the Yang--Baxter equation is a pair $(X,r)$,
where $X$ is a nonempty set and $r\colon X\times X\longrightarrow X\times X$
is a bijection satisfying
\[
 (r\times\id_X)(\id_X\times r)(r\times\id_X)
 =
 (\id_X\times r)(r\times\id_X)(\id_X\times r).
\]
Writing $r(x,y)=(\sigma_x(y),\tau_y(x))$, the solution is called
\emph{non-degenerate} if every $\sigma_x$ and every $\tau_y$ is a bijection,
and \emph{involutive} if $r^2=\id_{X\times X}$.
Its permutation group is
\[
 \mathcal G(X,r)=\langle\sigma_x\mid x\in X\rangle\leq\Sym(X).
\]
A finite group is an \emph{involutive Yang--Baxter group}, or an
\emph{IYB-group}, if it is isomorphic to the permutation group of a finite
involutive non-degenerate solution.

Left braces, introduced by Rump \cite{Rump2007}, provide an algebraic approach
to these groups. A \emph{left brace} is a set $B$ with two operations
$+$ and $\cdot$ such that $(B,+)$ is an abelian group, $(B,\cdot)$ is a
group, and
\[
 a(b+c)=ab+ac-a\qquad(a,b,c\in B).
\]
By \cite[Theorem 2.1]{CedoJespersRio2010}, a finite group is an IYB-group if and only if
it is isomorphic to the multiplicative group of a left brace.
Thus the IYB problem can be formulated as the problem of choosing a
compatible abelian group operation on a given group.

Finite IYB-groups are soluble \cite{EtingofSchedlerSoloviev1999}. The converse does not hold:
Bachiller \cite{Bachiller2016} constructed finite $p$-groups which are not
IYB-groups. On the positive side, finite nilpotent groups of class at
most two admit such structures; the corresponding circle-group
construction goes back to Ault and Watters \cite{AultWatters1973}.
Equivariant brace structures and their behaviour under group
factorisations were subsequently investigated in \cite{MengBallesterEstebanFuster2021,BEJP2023}.

Ced\'o and Okni\'nski \cite{CedoOkninski2025} proved that a finite soluble group whose
Sylow subgroups have nilpotency class at most two is an $\IYB$-group if its
Sylow $2$-subgroups are abelian. They asked whether the last assumption
can be omitted.
Meng, Ballester-Bolinches and Jim\'enez-Seral \cite{MengBallesterJimenez2026} established the conclusion when the nilpotent
residual is $Q_8$-free, or also when the Sylow $2$-subgroups are
isomorphic to the quaternion group $Q_8$.

In this paper we give an affirmative answer to Ced\'o and Okni\'nski's question without
any additional assumption.

\begin{mainthm}\label{thm:main}
Let $G$ be a finite soluble group whose Sylow subgroups have nilpotency
class at most two. Then $G$ is an $\IYB$-group.
\end{mainthm}

Recall that a group $A$ \emph{acts on a left brace} $(B,+,\cdot)$ if there is an action of $A$ on the set $B$ such that $a \cdot(g+h) = a\cdot g + a\cdot h$ and $a\cdot(gh) = (a\cdot g)(a\cdot h)$, for all $g, h \in B$, that is, if the action of $A$ on the set $B$ is actually an action of $A$ on the group $(B,+)$ and an action of $A$ on the group  $(B,\cdot)$.

An action of a group $A$ on an $\IYB$-group $(G, \cdot)$, for which it is understood that the associated left brace is $(G, {+}, {\cdot})$, is said to be \emph{equivariant} with respect to $\IYB$-structure $(G, {+}, {\cdot})$ if $A$ acts on the skew left brace $(G,+,\cdot)$. In this case, we say that $G$ has an \emph{$A$-equivariant} $\IYB$-structure.
For a left brace $B$, put
\[
\lambda_a(b)=-a+ab\qquad(a,b\in B).
\]
The maps $\lambda_a$ are additive automorphisms, and
$\lambda\colon(B,\cdot)\longrightarrow\Aut(B,+)$ is a homomorphism.
We write $\Ker(B)=\ker\lambda$.

It is well-known that, for an odd order nilpotent group $N$ of nilpotency class at most two, define
$$x+y:=xy[x,y]^{1/2},$$
where $[x,y]^{1/2}$ is the unique 
element in the odd order group $N$ satisfying the equation $u^2=[x,y]$.
Then the $\IYB$-structure $(N,+,\cdot)$ is $\Gamma$-equivariant for any action of a group $\Gamma$ on $N$ via automorphism. (e.g. See~\cite[Lemma~2.1]{CedoOkninski2025})
It is natural to see that the choice of $[x,y]^{1/2}$ will fail for even order groups. In fact, we have provided some other selection methods for $2$-groups.

\begin{mainthm}\label{thm:average}
Let a group $\Gamma$ act on a $2$-group $N$ of nilpotency class at most two by automorphisms. Write $K=\C_\Gamma(N/\Z(N))$ and suppose that $|\Gamma/K|$ is odd. Then $N$ has a $\Gamma$-equivariant $\IYB$-structure $(N,+,\cdot)$
such that $\Z(N)\leq\Ker(N)$.
\end{mainthm}

Note that, in Theorem~\ref{thm:average}, $K$ acts trivrally on $N/\Z(N)$. In some sense, $K$ acts on $N$ as central automorphisms. Hence Theorem~\ref{thm:average} is actually a generalized form of \cite[Proposition~2.5]{MengBallesterEstebanFuster2021}.

\emph{From now on, all groups considered will be finite.}

\section{Equivariant brace structures}

It is well-known that the use of symplectic bilinear forms is very effective for dealing with the structure and automorphisms of extraspecial $p$-groups, e.g. See \cite[Chapter~3, Theorem~13.7]{Huppert2025}. Here we need a slightly generalized bilinear form to handle general nilpotent groups of class at most two.

Let $V, C$ be two abelian groups. We say $\varphi \colon V\times V\longrightarrow C$ be a bilinear map if all
maps $\varphi(\cdot, v)$ and $\varphi(u, \cdot), u,v \in V$, are group-homomorphisms from $V$ to $C$, that is,
\[\varphi(u_1u_2,v)=\varphi(u_1,v)\varphi(u_2,v);\] \[\varphi(u,v_1v_2)=\varphi(u,v_1)\varphi(u,v_2)\]
holds for all $u_1,u_2,u,v_1,v_2,v \in V$.
Moreover, a bilinear map is called \emph{alternating} if
\(\varphi(u,u)=1\) for each $u \in V$.

Assume that a group $\Gamma$ act on $V$ and $C$ via automorphism respectively. A bilinear map $\varphi: V \times V \longrightarrow C$ is called
$\Gamma$-invariant if \[\gamma \cdot \varphi(u,v)=\varphi(\gamma \cdot u, \gamma \cdot v).\]
holds for all $\gamma \in \Gamma, u,v \in V$.

\begin{lemma}\label{lem:alternating}
Let a group $\Gamma$ act on two abelian groups $V$ and $C$ via automorphism respectively. Let
$c\colon V\times V\longrightarrow C$ be an $\Gamma$-invariant alternating bilinear map.
Write $K=\C_{\Gamma}(V) \cap \C_{\Gamma}(C)$ and suppose that
\[(|\Gamma/K|, \exp(C))=1.\]
Then there is a $\Gamma$-invariant bilinear map $\varphi\colon V\times V\longrightarrow C$
such that \(\varphi(u,v)\varphi(v,u)^{-1}=c(u,v) \) for all
\(u,v\in V.\)
\end{lemma}
\begin{proof}
Write $V=\langle e_1\rangle \times \langle e_2\rangle \cdots \times \langle e_t\rangle  $, where $o(e_i)=m_i$.
Define
\begin{equation}\label{eq:triangular}
 \varphi_0\left(\prod_i  e_i^{a_i},\prod_j  e_j^{b_j}\right)
 =\prod_{i<j}c(e_i,e_j)^{a_i b_j}.
\end{equation}
Since $c$ is bilinear, \(c(e_i,e_j)^{m_i}=c(e_i,e_j)^{m_j}=1.\)
Therefore the right-hand side of \eqref{eq:triangular} is independent
of the chosen integer $a_i, b_j$. Hence $\varphi_0$ is a well-defined map and it is easy to check that $\varphi_0$ is
bilinear. By hypothesis, $c$ is alternating, which implies that $c(e_j,e_i)=c(e_i,e_j)^{-1}$ and $c(e_i,e_i)=1$. It follows that
\begin{equation}\label{eq:alternating}
 c(u,v)=\prod_{i<j}c(e_i,e_j)^{a_i b_j-a_j b_i}=\varphi_0(u,v)\varphi_0(v,u)^{-1},
\end{equation}
where \(u=\prod_i  e_i^{a_i}\) and \(v=\prod_j  e_j^{b_j}\).

Write $m=|\Gamma/K|$ and and $\exp(C)=e$. Assume that $\Gamma/K=\{\gamma_1K, \cdots, \gamma_mK\}$, where $\{\gamma_1, \cdots \gamma_m\} \subseteq \Gamma$ is a set of left coset representatives of $K$ in $\Gamma$. As $(m,e)=1$, we may choose an integer $r$ such that
$rm\equiv1\pmod{e}$, and define
\begin{equation}\label{eq:average}
	\varphi(u,v)=
	\left(\prod_{1 \leq i \leq  m} \gamma_i^{-1} \varphi_0(\gamma_i u,\gamma_i v)\right)^r, u,v \in V.
\end{equation}
The values commute, and each factor is bilinear, so $\varphi$ is
bilinear. Moreover,
\[
\begin{aligned}
	\varphi(u,v)\varphi(v,u)^{-1}&=\left(\prod_{1 \leq i \leq  m} \gamma_i^{-1}\varphi_0(\gamma_i u,\gamma_iv)\right)^r \left(\prod_{1 \leq i \leq  m} \gamma_i^{-1}\varphi_0(\gamma_i v,\gamma_i u)\right)^{-r}\\
	&=\left(\prod_{1 \leq i \leq  m} \gamma_i^{-1} (\varphi_0(\gamma_i u,\gamma_i v) \varphi_0(\gamma_i v,\gamma_i u)^{-1})\right)^{r}\\
	&=\left(\prod_{1 \leq i \leq  m} \gamma_i^{-1}c(\gamma_i u, \gamma_i  v)\right)^{r}~~\text{(Applying Eq.\eqref{eq:alternating})}\\
	&=\left(\prod_{1 \leq i \leq  m}  c(u,v)\right)^{r} ~\text{(as $c$ is $\Gamma$-invariant)}\\
	&=c(u,v)^{rm}=c(u,v).
\end{aligned}
\]
Now we only have to show that $\varphi$ is $\Gamma$-invariant. For each $\gamma \in \Gamma$, as the right multiplication of $\gamma$ on $\Gamma/K$ permutates all cosets $K\gamma_1, \cdots, K \gamma_m$, we may assume that $K\gamma_i\gamma=K\gamma_{\alpha(i)}$, where $\alpha$ is a permutation on $\{1,2,\cdots, m\}$ (depending on $\gamma$).

As $K \unlhd \Gamma$, $\gamma_i\gamma K=\gamma_{\alpha(i)} K$ and set $\gamma_i\gamma=\gamma_{\alpha(i)} k_i$ for some $k_i \in K$. Hence it follows that
\[
\begin{aligned}
\varphi(\gamma u, \gamma v)&=	\left(\prod_{1 \leq i \leq  m} \gamma_i^{-1}\varphi_0(\gamma_i (\gamma u),\gamma_i (\gamma v))\right)^r\\
&=	\left(\prod_{1 \leq i \leq  m} \gamma_i^{-1}\gamma_{\alpha(i)}\gamma_{\alpha(i)}^{-1} \varphi_0((\gamma_i\gamma) u,(\gamma_i\gamma) v)\right)^r\\
&=	\left(\prod_{1 \leq i \leq  m} \gamma_i^{-1}\gamma_{\alpha(i)}k_i\gamma_{\alpha(i)}^{-1}\varphi_0((\gamma_{\alpha(i)}k_i) u,(\gamma_{\alpha(i)}k_i) v)\right)^r~\text{( $K$ acts trivially on $C, V$)}\\
&=	\left(\prod_{1 \leq i \leq  m} \gamma \gamma_{\alpha(i)}^{-1}\varphi_0(\gamma_{\alpha(i)} u,\gamma_{\alpha(i)}v)\right)^r=	\gamma\left(\prod_{1 \leq i \leq  m}  \gamma_{\alpha(i)}^{-1}\varphi_0(\gamma_{\alpha(i)} u,\gamma_{\alpha(i)}v)\right)^r\\
&=\gamma\varphi(u,v).
\end{aligned}
\]
The arbitrariness of $\gamma$ implies that $\varphi$ is $\Gamma$-invariant, as dersired.
\end{proof}

\begin{proof}[\textbf{\emph{Proof of Theorem B}}]
If $N$ is Abelian, take $x+ y=xy$ and the result is trivial. Now assume that $N$ is non-Abelian, put
$V=N/\Z(N)$, $C=N'$ and
\[
 c:V \times V \rightarrow C~\text{with}~c(\bar x,\bar y)=[x,y],
\]
\(\text{where}~\bar x=x\Z(N), \bar y=y\Z(N) \in N/\Z(N)~\text{and}~x,y \in N.\) It is easy to check by definition that \(c\) is a well-defined alternating bilinear map.

Since $\Z(N),N'$ are characteristic subgroups of $N$, $\Gamma$ can naturally act on $N'$ and $\Gamma$ also can act on $V=N/\Z(N)$ induced by the action of $\Gamma$ on $N$ as follows:
$$\gamma \cdot (x\Z(N)):= (\gamma \cdot x)\Z(N).$$
Hence, for each $k \in K=\C_{\Gamma}(N/\Z(N))$ and $x,y \in N$,
we may assume that $k x=xz_1$ and $k y=yz_1$ for some $z_1,z_2 \in \Z(G)$. It implies that
$k \cdot [x,y]=[kx,ky]=[xz_1,yz_2]=[x,y]$, and so $k$ fixes $N'=C$. Hence $K=\C_{\Gamma}(V) \cap \C_{\Gamma}(C)$. 
For each $\bar x, \bar y \in V$, it follows that
\[
	c(\gamma  \bar x, \gamma  \bar y)=c(\overline{\gamma  x}, \overline{\gamma  y})
=[\gamma x, \gamma  y]=\gamma  [x,y]=\gamma c(\bar x, \bar y).	
\]
Hence $c$ is $\Gamma$-invariant.
Since $|\Gamma/K|$ is odd and $C$ is a $2$-group, by Lemma~\ref{lem:alternating}, there is a $\Gamma$-invariant
bilinear map $\varphi \colon V\times V\longrightarrow C$
such that
\begin{equation*}
	\varphi(u,v)\varphi(v,u)^{-1}=c(u,v)
	\qquad(u,v\in V).
\end{equation*}
Now we define an addition "+" on $N$:
\[x+y:=xy\varphi(\bar x, \bar y)^{-1}\]
Since $xy=yx[x,y]$, it follows that
\[
\begin{aligned}
 y+x&=yx\,\varphi(\bar y, \bar x)^{-1}=yxc(\bar x, \bar y)\varphi(\bar x,\bar y)^{-1}\\
 &=yx[x,y]\varphi(\bar x,\bar y)^{-1}=xy\varphi(\bar x,\bar y)^{-1}\\
 &=x+y.
\end{aligned}
\]
Hence $+$ is commutative. Using bilinearity, we will see that both $(x+y)+z$ and
$x+(y+ z)$ are equal to
\[
 xyz\varphi(\bar x,\bar y)^{-1}\varphi(\bar x,\bar z)^{-1}\varphi(\bar y,\bar z)^{-1}.
\]
Thus $+$ is associative. For each $x \in N$, as $\varphi(\bar x, \bar 1)=\varphi(\bar x, \bar 1 \cdot \bar 1)=\varphi(\bar x, \bar 1)\varphi(\bar x, \bar 1)$, one obtain that $\varphi(\bar x, \bar 1)=1$. Hence $1+x=x+1=x 1 \varphi(\bar x, \bar 1)=x$ and $1$ is also the additive identity. Moreover, the additive
inverse of $x$ is $-x=x^{-1}\varphi(\bar x, \bar x)^{-1}$. Indeed, we have
\[
\begin{aligned}
x+(-x)&=x+x^{-1}\varphi(\bar x, \bar x)^{-1}\\
&=xx^{-1}\varphi(\bar x, \bar x)^{-1} \varphi(\bar x, \overline{x^{-1}\varphi(\bar x, \bar x)^{-1}})^{-1}\\
&=xx^{-1}\varphi(\bar x, \bar x)^{-1} \varphi(\bar x, \overline{x^{-1}})^{-1}~\text{(as $\varphi(\bar x, \bar x)^{-1} \in N' \leq \Z(N)$)}\\
&=\varphi(\bar x, \bar x)^{-1} \varphi(\bar x, \overline{x}^{-1})^{-1}=1.
\end{aligned}
\]
Therefore $(N,+)$ is an Abelian group. Define $\lambda_x(y)=(-x)+xy$ for $x,y \in N$. We see
$$
\begin{aligned}
\lambda_x(y)&=x^{-1}\varphi(\bar x, \bar x)^{-1}+xy\\
&=x^{-1}\varphi(\bar x, \bar x)^{-1}xy\varphi(\overline{x^{-1}}, \overline{xy})^{-1}\\
&=y\varphi(\bar x, \bar x)^{-1}\varphi(\overline{x}, \overline{xy})=y\varphi(\overline{x}, \overline{y}).
\end{aligned}$$
Then, for $y,z\in N$, we have
\begin{align*}
 \lambda_x(y+z)
 &=(y+z)\,\varphi(\bar x,\overline{y+z})=yz\,\varphi(\bar y, \bar z)^{-1}\varphi(\bar x,\bar y \bar z)\\
 &=y\varphi(\bar x,\bar y)z\varphi(\bar x,\bar z)\varphi(\bar y, \bar z)^{-1}\\
 &=y\varphi(\bar x,\bar y)z\varphi(\bar x,\bar z)\varphi(\overline{y\varphi(\bar x,\bar y)}, \overline{z\varphi(\bar x,\bar z)})^{-1}\\
 &=y\varphi(\bar x,\bar y)+z\varphi(\bar x,\bar z)=\lambda_x(y)+\lambda_x(z).
\end{align*}
It also implies that $x(y+z)=xy+xz-x$ for $x,y,z \in N$. Hence $(N,+,\cdot)$ is a left brace.
If $x\in\Z(N)$, then $\bar x=0$ and $\lambda_x(y)=y\varphi(\bar x,\bar y)=y$ for all $y \in N$.
Thus $\Z(N)\leq\Ker(N)$. Finally, the invariance of $\varphi$ gives
\[
 \gamma(x+y)=\gamma(xy \varphi(\bar x, \bar y)^{-1})
 =\gamma(x)\gamma(y)
   \varphi(\overline{\gamma(x)},\overline{\gamma(y)})^{-1}
 =\gamma(x)+\gamma(y),
\]
for every $\gamma\in\Gamma$. Hence $N$ has a $\Gamma$-equivariant $\IYB$-structure $(N,+,\cdot)$
such that $\Z(N)\leq\Ker(N)$ and the proof is complete.
\end{proof}

\begin{corollary}\label{cor:p-normalizer}
Let $G$ be a finite group whose Sylow $p$-subgroups have nilpotency
class at most two. If $N$ is a $p$-subgroup of $G$ and
$H\leq\N_G(N)$, then $N$ has an $H$-equivariant $\IYB$- structure
such that $\Z(N)\leq\Ker(N)$.
\end{corollary}

\begin{proof}
By \cite[Lemma~2.1]{CedoOkninski2025}, we only have to prove it when $p=2$. Take $S\in\Syl_2(H)$. Since $S$ normalises $N$, the product $NS$
is a $2$-subgroup of $G$, and therefore has class at most two. Hence
\[
 [N,S]\leq (NS)'\cap N\leq\Z(NS)\cap N\leq\Z(N),
\]
The kernel $K$ of the action of $H$ on
$N/\Z(N)$ contains $S$. Then $H/K$ is of odd order. Theorem~\ref{thm:average} applies.
\end{proof}

\begin{corollary}\label{cor:nilpotent}
Let $G$ be a finite group whose Sylow subgroups have nilpotency class
at most two. If $N$ is a nilpotent subgroup of $G$ normalised by
$H\leq G$, then $N$ has an $H$-equivariant $\IYB$-structure with $\Z(N)\leq\Ker(N)$.
\end{corollary}

\begin{proof}
Each Sylow subgroup $N_p$ of $N$ is characteristic in $N$ and is
normalised by $H$. Apply Corollary~\ref{cor:p-normalizer} and \cite[Lemma~2.1]{CedoOkninski2025} to $N_p$,
and take the direct product of the resulting braces. The action of
$H$ preserves each factor, and
\[
 \Z(N)=\prod_p\Z(N_p)\leq\prod_p\Ker(N_p)=\Ker(N),
\]
as desired.
\end{proof}

\section{Gluing brace structures}

We give the form of the gluing argument needed below. 
The first proposition is exactly the version of \cite[Theorem]{MengBallesterEstebanFuster2021} without the group action.
It is also the
left-brace specialisation of the factorisation results in \cite{BEJP2023}.

\begin{proposition}\label{prop:glue}
Let $G=NH$ be the product of $N$ and $H$ with $N\unlhd G$, and assume that $N,H$ have $\IYB$-structures $(N,+,\cdot)$ and $(H,+,\cdot)$ respectively. Suppose that
\begin{itemize}
	\item[(1)] $N\cap H \leq \Z(N)$ and $N\cap H \leq \Ker(N)\cap\Ker(H)$;
	\item[(2)] $(N,+,\cdot)$ is an $H$-equivariant $\IYB$-structure with respect to the conjugation action of $H$ on $N$.
\end{itemize}
Then $G$ has an $\IYB$-structure such that
\[\Ker(N)\,\C_{\Ker(H)}(N)\leq\Ker(G).\]
\end{proposition}

The following proposition is actually the
left-brace specialisation of \cite[Theorem C]{MengBallesterJimenez2026}. We give a proof here for completeness.

\begin{proposition}\label{prop:iterated-glue}
	Let $G=A_1\cdots A_k$ be the product of some subgroups $A_1,A_2,\cdots,A_k$, put $H_i=A_{i+1}\cdots A_k$ for $i<k$,
	and put $H_k=1$ and $A_{k+1}=1$. Suppose that each $A_i$ has an $\IYB$-structure $(A_i,+,\cdot)$ and that, for every $i \in \{1,2,\cdots, k\}$,
	\begin{enumerate}
		\item $H_i$ normalises $A_i$;
		\item $(A_i, +, \cdot)$ is an $H_i$-equivariant $\IYB$-structure with the respect to the action by conjugation of $H_i$ on $A_i$;
		\item $(A_1\cdots A_i)\cap H_i\leq\C_G(A_iA_{i+1})$;
		\item $\Z(A_i)\leq\Ker(A_i)$.
	\end{enumerate}
Then $G$ has an $\IYB$-structure $(G, +, \cdot)$ satisfying the following property:

(*) Suppose that  $a_i\in \Z(A_i)$ for each $i=1,2,\cdots k$, and $a_ia_{i+1}\cdots a_k\in \C_G(A_{i-1})$ for each $i=2, \cdots,k$. Then  $a_1a_2\cdots a_k \in \Ker(G)$.

Furthermore, $(A_i,+,\cdot)$ is a subbrace of $(G,+,\cdot)$ for every $i \in \{1, \cdots,k\}$.
\end{proposition}
\begin{proof}
	We argue by induction on $k$. It is clear that we may assume that  $k\geq 2$. It is rather clear that $H_1 = A_2\cdots A_{k}$ satisfies the hypotheses of the theorem. Then, by induction, $H_1$ has an $\IYB$-structure $(H_1,+,\cdot)$ satisfying Property~(*) and $(A_i,+,\cdot)$ are subbraces of $(H_1,+,\cdot)$ for every $i \in \{2, \cdots,k\}$. Furthermore, $G = A_1H_1$, $H_1$ normalises $A_1$ and the $\IYB$-structure
	$(A_1,+,\cdot )$ is $H_1$-equivariant with the respect to the action by conjugation of $H_1$ on
	$A_1$ by Statement~(2). Let  $a_1=a_2\cdots a_{k} \in A_1\cap H_1$, where
	$a_i \in A_i$,  $1\leq i \leq k$. Then $(a_1)^{-1}a_2\cdots a_k=1$. If $i \in \{2, \cdots, k\}$, then $(a_1)^{-1}a_2\cdots a_{i-1} a_i=(a_{i+1}\cdots a_{k})^{-1} \in A_1\cdots A_i\cap H_i\leq \C_G(A_iA_{i+1})\leq \C_G(A_{i})$ and $(a_1)^{-1}a_2\cdots a_{i-1} = (a_{i}\cdots a_{k})^{-1}\in \C_G(A_{i-1}A_{i})\leq \C_G(A_{i})$ by Statement~(3). Therefore, $a_{i}\in \Z(A_{i})$. Now, $(a_1)^{-1}a_2\cdots a_{i-1}=(a_ia_{i+1}\cdots a_{k})^{-1}\in \C_G(A_{i-1})$ and $a_ia_{i+1}\cdots a_{k}\in \C_{G}(A_{i-1})$. Since $H_1$ satisfies Property~(*), $a_2\cdots a_k\in \Ker (H_1)$. Note that  $A_1\cap H_1\leq \C_G(A_1A_2)\cap A_1\leq \Z(A_1)\leq \Ker (A_1)$ by Statement~(4). Consequently, $A_1 \cap H_1 \leq \Ker (A_1) \cap \Ker (H_1)$. Applying Proposition~\ref{prop:glue}, $G$ has an $\IYB$-structure $(G, +, \cdot)$ such that $\Ker(A_1)\C_{\Ker (H_1)}(A_1)\leq \Ker (G)$ and $(A_1, +, \cdot)$ and $(H_1,+,\cdot)$ are subbraces of $(G,+,\cdot)$. Hence $(A_i,+,\cdot)$ is a subbrace of $(G,+,\cdot)$, for every $i \in \{1, \cdots,k\}$.
	
	Assume that $a=a_1a_2\cdots a_k\in \Z(A_1)\Z(A_2)\cdots \Z(A_k)$, where $a_i\in \Z(A_i)$, and $a_ia_{i+1}\cdots a_k\in \C_G(A_{i-1})$, for each $i=1, \cdots,k$. Then $a_2\cdots a_k\in \C_G(A_1)$ and $a_2\cdots a_k\in \Ker (H_1)$. Thus $a_1a_2\cdots a_k\in \Ker(A_1)\C_{\Ker (H_1)}(A_1)\leq \Ker(G)$. This shows that $(G, +, \cdot)$ has Property~(*).	
\end{proof}

\section{Proof of Theorem~\ref{thm:main}}

We shall use the following standard properties in the precise form
established in \cite{MengBallesterJimenez2026}.
We first recall the definition of $\mathcal{N}$-decomposition in soluble groups \cite[Definition 11]{MengBallesterJimenez2026}.

\begin{definition}
	An \emph{$\mathcal{N}$-decomposition} of a group $G$ is a set of nilpotent subgroups $A_1,A_2,\cdots,A_k$, $k\geq 2$, of $G$ satisfiying the following properties:
	\begin{enumerate}
		\item $G=A_1A_2\cdots A_k$;
		\item $A_i$ is normalised by $A_j$ for each $1\leq i <j\leq k$.
		\item $(A_1\cdots A_i) \cap (A_{i+1} \cdots A_k) \leq A_1' \cdots A_i'$ for each $1 \leq i \leq k-1$.
	\end{enumerate}
\end{definition}

This is a partial result of \cite[Theorem~12]{MengBallesterJimenez2026}.
\begin{theorem}\label{Ndec}
	Assume that the action of a group  $A$ on a group $G$ is coprime. Then $G$ is soluble if and only if $G$ has an $\mathcal{N}$-decomposition $\{A_1,A_2,\cdots,A_k\}$  composed of $A$-invariant subgroups $A_i$, $1 \leq i \leq k$. 
\end{theorem}

The following result is a special case of \cite[Lemma~15]{MengBallesterJimenez2026}.
\begin{lemma}\label{lem-nilpotency-class-2}
Let $G$ be a soluble group such that all Sylow subgroups of $G$ have nilpotency class at most two. Let $\{A_1,A_2,\cdots,A_k\}$ be an $\mathcal{N}$-decomposition of $G$. Then
	$$A_1A_2\cdots A_i \cap A_{i+1}\cdots A_k \leq \C_G(A_iA_{i+1}) \text { for each } 1 \leq i \leq k-1.$$
\end{lemma}

Now we can prove the main theorem.

\begin{proof}[\textbf{\emph{Proof of Theorem~\ref{thm:main}}}]
Since $G$ is soluble and all Sylow subgroups of $G$ have nilpotency class at most two, by Theorem~\ref{Ndec} and Lemma~\ref{lem-nilpotency-class-2}, $G$ has an $\mathcal{N}$-decomposition $\{A_1,A_2,\cdots,A_k\}$ such that
$A_1A_2\cdots A_i \cap A_{i+1}\cdots A_k \leq \C_G(A_iA_{i+1})$ for each $1\leq i \leq k-1$.
Write $H_i=A_{i+1}A_{i+2} \cdots A_k$ for each $1 \leq i \leq k-1$ and $H_k=1$.
For each $i$, 
$A_i$ is nilpotent and is normalised by $H_i$.
Corollary~\ref{cor:nilpotent} gives an $\IYB$-structure on $A_i$
which is $H_i$-equivariant and satisfies \( \Z(A_i)\leq\Ker(A_i).\)
All four hypotheses of
Proposition~\ref{prop:iterated-glue} are  satisfied, and hence $G$ is an $\IYB$-group.
\end{proof}

\textbf{Acknowledgment:} The first author is sponsored by Natural Science Foundation of Shanghai
(24ZR1422800) and National Natural Science Foundation of China (12471018). The second author is supported by the grant PID2024-159495NB-I00, from MICIU/AEI/10.13039/501100011033 and ERDF/EU. He is also supported by the grant CIAICO/2023/007 from the Conselleria d'Educaci\'o, Universitats i Ocupaci\'o, Generalitat Valenciana and by the program "High-end Forcing Expert Program of China(H20240848)", Shanghai University.

\bibliographystyle{plain}
\bibliography{bibM}

\begin{thebibliography}{10}

\bibitem{AultWatters1973}
J.~C. Ault and J.~F. Watters.
\newblock Circle groups of nilpotent rings.
\newblock {\em Amer. Math. Monthly}, 80:48--52, 1973.

\bibitem{Bachiller2016}
D.~Bachiller.
\newblock Counterexample to a conjecture about braces.
\newblock {\em J. Algebra}, 453:160--176, 2016.

\bibitem{BEJP2023}
A.~Ballester-Bolinches, R.~Esteban-Romero, P.~Jim\'{e}nez-Seral, and
  V.~P\'{e}rez-Calabuig.
\newblock On {Y}ang-{B}axter groups.
\newblock {\em Quaest. Math.}, 46(7):1273--1281, 2023.

\bibitem{CedoJespersRio2010}
F.~Ced{\'o}, E.~Jespers, and A.~Del R{\'\i}o.
\newblock Involutive {Y}ang-{B}axter groups.
\newblock {\em Trans. Amer. Math. Soc.}, 362(5):2541--2558, 2010.

\bibitem{CedoOkninski2025}
F.~Ced\'{o} and J.~Okni\'{n}ski.
\newblock New classes of {IYB} groups.
\newblock ArXiv:2502.14805, 2025.

\bibitem{EtingofSchedlerSoloviev1999}
P.~Etingof, T.~Schedler, and A.~Soloviev.
\newblock Set-theoretical solutions to the quantum {Y}ang-{B}axter equation.
\newblock {\em Duke Math. J.}, 100(2):169--209, 1999.

\bibitem{Huppert2025}
B.~Huppert.
\newblock {\em Finite groups I}.
\newblock Springer, 2025.
\newblock Translated by Christopher A. Schroeder.

\bibitem{MengBallesterEstebanFuster2021}
H.~Meng, A.~Ballester-Bolinches, R.~Esteban-Romero, and N.~Fuster-Corral.
\newblock On finite involutive {Y}ang-{B}axter groups.
\newblock {\em Proc. Amer. Math. Soc.}, 149(2):793--804, 2021.

\bibitem{MengBallesterJimenez2026}
H.~Meng, A.~Ballester-Bolinches, and P.~Jim{\'e}nez-Seral.
\newblock On a class of involutive {Yang--Baxter} groups.
\newblock {\em J. Algebra}, 704:1--16, 2026.

\bibitem{Rump2007}
W.~Rump.
\newblock Braces, radical rings, and the quantum {Y}ang-{B}axter equation.
\newblock {\em J. Algebra}, 307:153--170, 2007.

\end{thebibliography}
\end{document}